\documentclass[12pt]{amsart}
\usepackage{amssymb, amscd, latexsym, color, amscd, float}
\usepackage{url}
\usepackage[final]{graphicx}

\usepackage{ulem,ifthen}
\newcommand{\ifemptythenelse}[3]{%
  \begingroup
    \def\dummy{#1}%
    \def\empty{}%
    \ifx\dummy\empty{#2}\else{#3}\fi
  \endgroup
  }
\DeclareRobustCommand{\change}[2][]{%
  \ifemptythenelse{#1}{%
    \ifemptythenelse{#2}{}{\begin{color}{blue}{{#2}}\end{color}}%+
  }{%
    \ifemptythenelse{#2}{%
      \sout{\begin{color}{red}{{#1}}\end{color}}%
    }{%
      \ifmmode
        \begin{color}{green}{{#2}}\end{color}%
      \else
      \begin{color}{blue}{{#2}}\end{color}%
      \footnote{was: \begin{color}{red}{{#1}}\end{color}}%
      \fi
    }%
  }%
}

\newcommand{\washere}[1]{}

\newtheorem{thm}{Theorem}[section]
\newtheorem{prop}[thm]{Proposition}
\newtheorem{lem}[thm]{Lemma}
\newtheorem{cor}[thm]{Corollary}

\newtheorem{defn}[thm]{Definition}
\newtheorem{rmk}[thm]{Remark}

\begin{document}

\title
[]
%[\normalsize A generalization of periodicity theorem]
{\normalsize
On finite systems of equations\\ in acylindrically hyperbolic groups
%Verbally closed acylindrically hyperbolic subgroups are algebraially closed
}
\bigskip

\author{Oleg Bogopolski}
\address{{Sobolev Institute of Mathematics of Siberian Branch of Russian Academy
of Sciences, Novosibirsk, Russia}\newline {and D\"{u}sseldorf University, Germany}}
\email{Oleg$\_$Bogopolski@yahoo.com}

\begin{abstract}
%Oleg Marie Rafael Vasilii
Let $H$ be an acylindrically hyperbolic group without nontrivial finite normal subgroups.
We show that any finite system $S$ of equations with constants from $H$ is equivalent to a single equation.

%An equation of the form
We also show that the algebraic set associated with $S$ is, up to conjugacy, a projection of
the algebraic set associated with a single splitted equation (such equation has the form $w(x_1,\dots,x_n)=h$, where $w\in F(X)$, $h\in H$).

From this we deduce the following statement:
{\it Let $G$ be an arbitrary overgroup of the above group $H$. Then
$H$ is verbally closed in $G$ if and~only~if it is algebraically closed in $G$.}

Another corollary: If $H$ is a non-cyclic torsion-free hyperbolic group,
then every (possibly infinite) system of equations with finitely many variables and with constants from $H$
is equivalent to a single equation.
\end{abstract}

\maketitle

\setcounter{tocdepth}{1}
%\tableofcontents

\bigskip

%\begin{center}
%{\bf \Large Equations in acylindrically hyperbolic groups
%and verbal closedness}

%\medskip

%(Oleg Bogopolski)
%\end{center}

%\bigskip

\section{Introduction}

Study of equations in groups is a classical subject of group theory with
long history, nice results and many intriguing problems~(see the inspiring paper of Neumann~\cite{Neumann_0} of 1943, the survey of
Roman'kov~\cite{Roman'kov} of 2012 and the recent paper of Klyachko and Thom~\cite{KT}).

Acylindrically hyperbolic groups (definied by Osin in~\cite{Osin_1}, see Definition~\ref{Definition_of_Osin} below) is a large class of
groups that is intensively studied in the modern geometric group theory.
This class includes non-(virtually cyclic) groups
%~\cite{Gromov}
that are hyperbolic relative to proper subgroups,
%~\cite{Osin_0},
many 3-manifold groups, groups of deficiency at least~2, many groups acting on trees,
non-(virtually cyclic) groups acting properly on proper CAT(0)-spaces and containing rank-one elements, non-cyclic directly indecomposable right-angled Artin groups,
all but finitely many mapping class groups,
% of punctured closed surfaces,
$Out(F_n)$ for $n\geqslant 2$,
and many other interesting groups; see the survey of Osin~\cite{Osin_3}.

However almost nothing is known on solutions of equations and related problems
in the class of acylinfrically hyperbolic groups.
In~\cite{Bog_1}, we described solutions of certain equations of the form $x^ny^m=a^nb^m$ in acylindrically hyperbolic groups. Using this description, we studied in~\cite{Bog_1} the {\it verbal closedness}  (see Definition~\ref{alg,verb,retr}  below) of acylindrically hyperbolic subgroups in groups.

\newpage

One of the corollaries there
solves Problem 5.2 from the paper~\cite{MR} of Myasnikov and Roman'kov:
{\it Verbally closed subgroups of torsion-free hyperbolic groups are retracts.}

More information on algebraic and verbal closedness of groups in overgroups or in classes of groups can be found in~\cite[Sections 1, 2, 15]{Bog_1}.
The present paper is a continuation of~\cite{Bog_1}.
We formulate briefly the main results.
For convenience, we call a group {\it clean} if it does not contain nontrivial finite normal subgroups.

\medskip

$\bullet$ Theorem {\bf A}, says that {\it if $H$ is a clean acylindrically hyperbolic group, then
any finite system of equations with constants in $H$ has the same set
of solutions in $H$ as a single equation}. Moreover, this set is a projection (up to conjugacy) of the set of solutions of a single {\it splitted} equation (see Definition~\ref{def_splitted} below).

\medskip

$\bullet$ Theorem {\bf C}, says that
{\it for any clean acylindrically hyperbolic group $H$ and any overgroup $G$ of $H $
the notions of verbal and algebraic closedness of $H$ in~$G$ are equivalent.}
%The proof of this theorem is short modulo Theorem~A.
Some special cases of this theorem
%(on virtially free groups and on free products of groups)
were considered earlier in~\cite{KMM} and~\cite{Mazhuga_3}, see Remark~\ref{Remark_free_prod} below.

\medskip

$\bullet$ A part of Corollary~{\bf D} says that if
$H$ is a finitely generated clean acylindrically hyperbolic group and $G$ is a finitely presented overgroup of $H$,
%additionally, $H$ is finitely generated and $G$ is finitely presented,
then the notions for $H$ to be verbally closed in $G$, to be algebraically closed in $G$, and to be a retract
of $G$ are equivalent.

\medskip

$\bullet$ A special case of Corollary {\bf B} says that if $H$ is a non-cyclic torsion-free hyperbolic group,
then every (possibly infinite) system of equations with constants from $H$ and finitely many variables
is equivalent to a single equation with constants from $H$, i.e., they have the same sets of solutions in~$H$.

\medskip

In the proof we use test words for acylindrically hyperbolic groups, see~\cite{Bog_1} and Section~5 below.

%Zarisski topology on $G^{\mathbb{Z}}$.

\section{Full formulations of main results and definitions}

\medskip

%Let $F_n$ be the free group with basis $x_1,\dots,x_n$.
Let $H$ be a group. An {\it equation} with variables $x_1,\dots,x_n$ and constants from~$H$ is an
element of the free product
$F_n\ast H$, where $F_n$ is the free group  with basis $x_1,\dots,x_n$.
Sometimes we write an equation $f$ in the form $f(x_1,\dots,x_n;H)$ stressing that $f$ involves the variables $x_1,\dots,x_n$ and constants from $H$.
Sometimes, for convenience, we write an equation $f_0f_1$ in the form $f_0=f_1^{-1}$.

\medskip

%If the left side of the equation does not contain letters from $H$, we will omit~$H$ from this expression.
Let $S\subseteq F_n\ast H$ be a system of equations and let $G$ be an overgroup of $H$.\break
A tuple $(g_1,\dots ,g_n)$ with components from $G$
is called a {\it solution of the system $S$ in} $G$ if $f(g_1,\dots,g_n;H)=1$ in $G$ for every
equation $f(x_1,\dots,x_n;H)$ from~$S$.
Let $V_G(S)$ be the set of all solutions of the system $S$ in $G$, i.e.,
$$
V_G(S)=\{(g_1,\dots ,g_n)\in G^n\,|\,  f(g_1,\dots,g_n;H)=1\hspace*{1mm}{\text{ \rm for all}}\hspace*{1mm} f\in S\}.
$$

%\newpage

\subsection{Systems of equations versus a single equation}

%Let $F_n$ be the free group of rank $n$ with basis $x_1,\dots,x_n$. Let $F(X)$ be the free group with basis $X$.

\begin{defn}\label{def_splitted}
{\rm
An equation $f\in F_n\ast G$ is called {\it splitted} if it has the form $wg$, where $w\in F_n$ and $g\in G$.
}
\end{defn}

%For a system of equations $S\subseteq F_n\ast G$, let $V_G(S)$ be the set of all solutions
%of the system $S$ in $G$, i.e.
%$$
%V_G(S)=\{(g_1,\dots ,g_n)\in G^n\,|\,  f(g_1,\dots,g_n)=1\hspace*{1mm}{\text{ \rm for all}}\hspace*{1mm} f\in S\}.
%$$

For $m\geqslant n$, let ${\text{\bf pr}}_n:G^m\rightarrow G^n$ be the projection to the first $n$ coordinates,
i.e., ${\text{\bf pr}}_n(g_1,\dots,g_m)=(g_1,\dots,g_n )$. For $g,u\in G$, we denote $g^u=u^{-1}gu$.\break
For $(g_1,\dots,g_n)\in G^n$ and $u\in G$, we set
$(g_1,\dots,g_n)^u=(g_1^u,\dots,g_n^u)$.
The first main result of this paper is the following theorem.

\medskip

\noindent
{\bf Theorem A.}
%\begin{thm}\label{algebraic_sets}
{\it
Let $H$ be an acylindrically hyperbolic group without nontrivial finite normal subgroups.
Let $S\subset F_n\ast H$ be a finite system of equations with constants from $H$.
Then the following statements hold:

\begin{enumerate}
\item[(1)] There exists a single equation $f\in F_n\ast H$ such that
$$
V_H(f)=V_H(S).
$$

\medskip

\item[(2)] There exists a natural $k\geqslant n$ and a single {\rm splitted} equation
$f\in F_k\ast H$ of the form $f_1f_0$, where $f_1\in F_k$ and $f_0\in H$ such that the following two properties are satisfied:

\begin{enumerate}
\item[(a)]
$$
{\text{\bf pr}}_n\bigl(V_H(f)\bigr)=\underset{\alpha\in \mathbb{Z}}{\bigcup} V_H(S)^{f_0^{\alpha}}.
$$

\item[(b)] For any overgroup $G$ of the group $H$ we have
$$
{\text{\bf pr}}_n\bigl(V_G(f)\bigr)\supseteq \underset{\alpha\in\mathbb{Z}}{\bigcup}V_G(S)^{f_0^{\alpha}}.
$$
\end{enumerate}
\item[(3)] There exist a natural $k\geqslant n$ and two splitted equations
$f,g\in F_k\ast H$ such that
$$
V_H(S)={\text{\bf pr}}_n\bigl(V_H(f)\bigr)\cap {\text{\bf pr}}_n\bigl(V_H(g)\bigr).
$$
\end{enumerate}
}
%\end{thm}

\medskip

Note that in the proof of this theorem we essentially use a proposition on test words in acylindrically hyperbolic groups obtained in~\cite{Bog_1}.

\medskip

Recall that a group $H$ is called {\it equationally noetherian} if
every system of equations with constants from $H$ and a finite number of variables
is equivalent to a finite subsystem, see~\cite{BMR}.
Since hyperbolic groups are equationally noetherian (see~\cite[Corolary 6.13]{RW} for the general case and
~\cite[Theorem 1.22]{Sela_2} for the torsion-free case),
the following corollary follows directly from Theorem~A.

% The torsion-free case was considered by Sela~\cite[Theorem 1.22]{Sela_2} and the general case by Reinfeldt and Weidmann~\cite[Corolary 6.13]{RW} (they expanded Sela's methods).
%From this and statement Theorem A follows that if $H$ is hyperbolic group without nontrivial finite normal
%subgroups, then every system of equations with constants from $H$ is equivalent to a single equation.
%More precisely:

\medskip

\noindent
{\bf Corollary B.}
%\begin{thm}\label{algebraic_sets}
%{\it
%Let $H$ be a hyperbolic group without nontrivial finite normal subgroups.
%Let $S\subseteq F_n\ast H$ be a (possibly infinite) system of equations with constants from $H$.
%Then there exists a single equation $f\in F_n\ast H$ such that
%$
%V_H(f)=V_H(S).
%$
%}
{\it
Let $H$ be a non-(virtually cyclic) hyperbolic group without nontrivial finite normal subgroups.
Then every (possibly infinite) system of equations with constants from $H$ and finitely many variables
is equivalent to a single equation with constants from $H$, i.e., they have the same sets of solutions in~$H$.
}

\medskip
%\bigskip

\subsection{Verbal and algebraic closedness.}

Let $X=\{x_1,x_2,\dots \}$ be a countably infinite set of variables and let $F(X)$ be the free group with basis $X$.\break 
We recall definitions of algebraically (verbally) closed subgroups and retracts.

\begin{defn}\label{alg,verb,retr} {\rm Let $H$ be a subgroup of a group $G$.

\begin{enumerate}
\item[(a)] (see~\cite{Neumann_2,MR})
The subgroup $H$ is called {\it algebraically closed} in $G$ if
for any finite system of equations
$$
S=\{W_i(x_1,\dots,x_n;H)=1\,|\, i=1,\dots,m\}
$$
with constants from $H$ the following holds: if $S$ has a solution in $G$, then it has a solution in $H$.

\item[(b)] (see~\cite[Definition 1.1]{MR}) The subgroup $H$ is called {\it verbally closed} in $G$ if
for any word $W\in F(X)$ and any element $h\in H$ the following holds: if the equation $W(x_1,\dots,x_n)=h$
has a solution in $G$, then it has a solution in $H$.

\item[(c)] The subgroup $H$ is called a {\it retract} of $G$ if there is a homomorphism
$\varphi:G\rightarrow H$ such that $\varphi|_{H}={\rm id}$. The homomorphism $\varphi$ is called a {\it retraction}.
\end{enumerate}
}
\end{defn}

%\medskip

Obviously, if $H$ is a retract of $G$, then $H$ is algebraically closed in $G$.
%The converse is true if $H$ is finitely generated and $G$ is finitely presented
%(see~\cite[Proposition 2.2]{MR}).
Algebraic closedness implies verbal closedness, but the converse implication is not valid in general, see
example in~\cite[Remark 13.2]{Bog_1}.
%An example of a verbally closed but not algebraically closed finitely generated subgroup
%of a finitely generated virtually free group is given in Remark~\ref{example}.

%\medskip

The following proposition of Myasnikov and Roman'kov says that, under some general assumptions,
the property of $H$ to be algebraically closed in $G$
is equivalent to the property of $H$ to be a retract of $G$.

%\medskip

Recall that a group $G$ is called {\it finitely generated over} a subgroup $H$ if there exists a finite subset $X\subseteq G$ such that $G=\langle X,H\rangle$.

\begin{prop}\label{alg_closed_retracts} {\rm (\cite[Proposition 2.2]{MR})}
Let $H$ be a subgroup of a group $G$. Suppose that at least one of the following holds:
\begin{enumerate}
\item[(a)] $H$ is finitely generated and $G$ is finitely presented.

\item[(b)] $H$ is equationally noetherian and $G$ is finitely generated over $H$.
\end{enumerate}
Then $H$ is algebraically closed in $G$ if and only if $H$ is a retract of $G$.
\end{prop}

%\medskip

Our second main theorem establishes the equivalence of verbal and algebraic closedness for
clean acylindrically hyperbolic subgroups of arbitrary groups.

\medskip

\noindent
{\bf Theorem C.}
%\begin{thm}\label{verbally_closed_implies_alg_closed}
{\it
Let $H$ be an acylindrically hyperbolic group without nontrivial finite normal subgroups
and let $G$ be an arbitrary overgroup of $H$.
Then $H$ is verbally closed in $G$ if and only if $H$ is algebraically closed in $G$.
}
%\end{thm}

\medskip

The proof of this theorem is very short modulo the second statement of Theorem {\bf A}, see Section~7.
The following corollary follows directly from Theorem~{\bf C} and Proposition~\ref{alg_closed_retracts}.
This corollary was already obtained by the author in~\cite{Bog_1} without appealing to Theorem {\bf C}.

\medskip

\noindent
{\bf Corollary D.}\label{verb_closed_retract_0} {\rm (\cite[Theorems 2.2 and 2.4]{Bog_1})}
Let $H$ be a subgroup of a group $G$ such that at least one of the following holds:
\begin{enumerate}
\item[(a)] $H$ is finitely generated and $G$ is finitely presented.

\item[(b)] $H$ is equationally noetherian and $G$ is finitely generated over $H$.
\end{enumerate}

Suppose additionally that $H$ is acylindrically hyperbolic
and does not have nontrivial finite normal subgroups.
Then the following three statements are equivalent.

\begin{enumerate}
\item[(1)] $H$ is algebraically closed in $G$.

\item[(2)] $H$ is verbally closed in $G$.

\item[(3)] $H$ is a retract of $G$.
\end{enumerate}
%\end{cor}

\section{Necessary definitions}

\subsection{Definitions of acylindrically hyperbolic groups}

\medskip

All actions of \break groups on metric spaces are assumed to be isometric in this paper.

\begin{defn} {\rm (see~\cite{Bowditch} and Introduction in~\cite{Osin_1})
An action of a group $G$ on a metric space $S$ is called
{\it acylindrical}
if for every $\varepsilon>0$ there exist $R,N>0$ such that for every two points $x,y$ with $d(x,y)\geqslant R$,
there are at most $N$ elements $g\in G$ satisfying
$$
d(x,gx)\leqslant \varepsilon\hspace*{2mm}{\text{\rm and}}\hspace*{2mm} d(y,gy)\leqslant \varepsilon.
$$
}
\end{defn}

Given a generating set $X$ of a group $G$, we say that the right Cayley graph $\Gamma(G,X)$ is
{\it acylindrical} if the left action of $G$ on $\Gamma(G,X)$ is acylindrical.
For Cayley graphs, the acylindricity condition can be rewritten as follows:
for every $\varepsilon>0$ there exist $R,N>0$ such that for any $g\in G$ of length $|g|_X\geqslant R$
we have
$$
\bigl|\{f\in G\,|\, |f|_X\leqslant \varepsilon,\hspace*{2mm} |g^{-1}fg|_X\leqslant \varepsilon \}\bigr|\leqslant N.
$$

Recall that an action of a group $G$ on a hyperbolic space $S$ is called {\it elementary} if the limit set
of $G$ on the Gromov boundary $\partial S$ contains at most 2 points.

\begin{defn}\label{Definition_of_Osin} {\rm (see~\cite[Definition 1.3]{Osin_1})
A group $G$ is called {\it acylindrically hyperbolic} if it satisfies one of the following equivalent
conditions:

\begin{enumerate}
\item[(${\rm AH}_1$)] There exists a generating set $X$ of $G$ such that the corresponding Cayley graph $\Gamma(G,X)$
is hyperbolic, $|\partial \Gamma (G,X)|>2$, and the natural action of $G$ on $\Gamma(G,X)$ is acylindrical.

\medskip

\item[(${\rm AH}_2$)] $G$ admits a non-elementary acylindrical action on a hyperbolic space.
\end{enumerate}
}
\end{defn}

In the case (AH$_1$), we also write that $G$ is {\it acylindrically hyperbolic with respect to $X$}.

\medskip

%Another equivalent definition of an acylindrically hyperbolic group uses
%the notion ``hyperbolically embedded subgroup'' (see ...).
%This notion appears in the cited below Lemmas~.... However we decided not to recall definitions of

\subsection{Hyperbolically embedded subgroups}
We use this notion only in formulations of lemmas in Section 4.
For convenience of the reader, we recall relevant definitions from~\cite{DOG}; see also~\cite{Osin_1}.

Let $G$ be a group with a fixed collection of subgroups $\{H_{\lambda}\}_{\lambda\in \Lambda}$.
Given a symmetric subset $X\subseteq G$ such that $G$ is generated by $X$ together with the union of all $H_{\lambda}$,
we denote by $\Gamma(G,X\sqcup \mathcal{H})$ the right Cayley graph of $G$ whose edges are labelled by letters from the alphabet $X\sqcup \mathcal{H}$, where
$$
\mathcal{H}=\bigsqcup_{\lambda\in\Lambda}H_{\lambda}.
$$
We consider the Cayley graph
$\Gamma(H_{\lambda},H_{\lambda})$ as a complete subgraph of $\Gamma(G,X\sqcup~\mathcal{H})$.

\begin{defn}{\rm (see~\cite[Definition 4.2]{DOG})}
{\rm For every $\lambda\in \Lambda$, we introduce a {\it relative metric}
$\widehat{d}_{\lambda}:H_{\lambda}\times H_{\lambda}\rightarrow [0,+\infty]$ as follows.
Let $a,b\in H_{\lambda}$. A path
in $\Gamma(G,X\sqcup \mathcal{H})$ from $a$ to $b$ is called {\it $H_{\lambda}$-admissible} if it has no edges in the subgraph $\Gamma(H_{\lambda},H_{\lambda})$.
The distance $\widehat{d}_{\lambda}(a,b)$ is defined to be the length of a shortest
{\it $H_{\lambda}$-admissible} path connecting $a$ to $b$ if such exists.
If no such path exists, we set $\widehat{d}_{\lambda}(a,b)=\!
\infty$.
}
\end{defn}

\begin{defn}\label{def_hyperb_embedd} {\rm (see~\cite[Definition 4.25]{DOG})
Let $G$ be a group, $X$ a symmetric subset of $G$. A collection of subgroups $\{H_{\lambda}\}_{\lambda\in \Lambda}$
of $G$ is called {\it hyperbolically embedded in $G$ with respect to $X$}
(we write $\{H_{\lambda}\}_{\lambda\in \Lambda}\hookrightarrow_h (G,X)$) if the following hold.

\begin{enumerate}
\item[(a)] The group $G$ is generated by $X$ together with the union of all $H_{\lambda}$ and the Cayley graph
$\Gamma(G,X\sqcup \mathcal{H})$ is hyperbolic.

\item[(b)] For every $\lambda\in \Lambda$, the metric space $(H_{\lambda},\widehat{d}_{\lambda})$ is
proper. That is, any ball of finite radius in $H_{\lambda}$ contains finitely many elements.
\end{enumerate}

\medskip

Further, we say that $\{H_{\lambda}\}_{\lambda\in \Lambda}$ is {\it hyperbolically embedded}
in $G$
%and write $\{H_{\lambda}\}_{\lambda\in \Lambda}\hookrightarrow_h G$
if $\{H_{\lambda}\}_{\lambda\in \Lambda}\hookrightarrow_h (G,X)$ for some $X\subseteq G$.
}
\end{defn}

It was proved in~\cite[Theorem 1.2]{Osin_1} that a group $G$ is acylindrically hyperbolic if and only if
it contains a proper infinite hyperbolically embedded subgroup.

\subsection{Elliptic and loxodromic elements}
The following definition is standard.
\begin{defn}
{\rm
Given a group $G$ acting on a metric space $S$, an element $g\in G$ is called {\it elliptic}
if some (equivalently, any) orbit of $g$ is bounded, and {\it loxodromic} if the map
$\mathbb{Z}\rightarrow S$ defined by
$n\mapsto g^nx$ is a quasi-isometric embedding for some (equivalently, any) $x\in S$. That is,
for $x\in S$, there exist $\varkappa\geqslant 1$ and $\varepsilon\geqslant 0$ such that for any $n,m\in \mathbb{Z}$ we have
$$
d(g^nx,g^mx)\geqslant \frac{1}{\varkappa} |n-m|-\varepsilon.
$$

Let $X$ be a generating set of $G$.
We say that $g\in G$ is {\it elliptic ({\rm respectively}, loxodromic) with respect to $X$} if $g$ is elliptic (respectively, loxodromic) for the canonical left action of $G$ on the Cayley graph $\Gamma(G,X)$.
%If $X$ is clear from a context, we omit the words ``with respect to $X$''.

%The set of all elliptic
%(respectively loxodromic) elements of $G$ with respect to $X$
%is denoted by ${\rm Ell}(G,X)$ (respectively by ${\rm Lox}(G,X))$.
}
%denote the sets of elliptic and loxodromic elements for the left action of $G$ on $\Gamma(G,X)$, respectively.}
\end{defn}

Note that even in the case of groups acting on hyperbolic spaces, there may be other types of actions
(see~\cite[Section 8.2]{Gromov} and~\cite[Section 3]{Osin_1}).
However, if $G$ is a group acting acylindrically on a hyperbolic space,
then every element of $G$ is either elliptic or loxodromic. This was first proved Bowditch~\cite[Lemma 2.2]{Bowditch}; a more general statement was proved by Osin in~\cite[Theorem 1.1]{Osin_1}).

\medskip

Recall that if $G$ is an acylindrically hyperbolic group with respect to a generating set $X$,
then every loxodromic, with respect to $X$, element $g\in G$ is contained in a
%any loxodromic element $g$ in an acylindrically hyperbolic group $G$ is contained in a
unique maximal virtually cyclic subgroup $E_G(g)$ of $G$, see\break \cite[Lemma 6.5]{DOG}. This subgroup is called the {\it elementary subgroup associated with~$g$}; it can be described as follows (see equivalent definitions in~\cite[Corollary~6.6]{DOG}):
$$
E_G(g)=\{f\in G\,|\, \exists  n\in \mathbb{N}:  f^{-1}g^nf=g^{\pm n}\}.\eqno{(3.1)}
$$
Clearly, the centralizer $C_G(g)$ of $g$ is contained in $E_G(g)$.

Recall that two elements $a,b\in G$ of infinite order are called {\it commensurable}
if there exist $g\in G$ and $s,t\in \mathbb{Z}\setminus \{0\}$ such that $a^s=g^{-1}b^tg$.

Preparing to the proof of Theorem A, we need to produce many non-commensurable loxodromic elements $g$ with respect to the {\it same} generating set of $G$ and with the additional property that $E_G(g)=\langle g\rangle$. Moreover, we need that these elements have certain form.
This technical part of the proof is the subject of Section~4.

\medskip

\section{Special elements}

We reproduce the following definition from~\cite{Bog_1}.

\begin{defn}\label{exact_defn_special}
{\rm  Suppose that $G$ is an acylindrically hyperbolic group.

\begin{enumerate}
\item[(a)] An element $g\in G$ is called {\it special} if there exists a generating set $X$ of $G$
such that

- $G$ is acylindrically hyperbolic with respect to $X$,

- $g$ is loxodromic with respect to $X$, and

- $E_G(g)=\langle g\rangle$.

\noindent
In this case $g$ is called {\it special with respect to~$X$}.

\medskip

\item[(b)] Elements $g_1,\dots ,g_k\in G$ are called {\it jointly special} if there exists a generating set $X$ of $G$
such that each $g_i$ is special with respect to $X$.
\end{enumerate}
}

\end{defn}

\medskip

Note that point (a) of this definition was already used in the case of relatively hyperbolic groups
(see comments in~\cite[Section 3]{OT}).

The purpose of this section is to prove Proposition~\ref{many_jointly_special}; it will be used
in Section 5 to construct certain test words in acylindrically hyperbolic groups.
%This proposition produces many non-commensurable jointly special
%elements in acylindrically hyperbolic groups. In Section ... we construct test words from these elements.
We deduce this proposition formally from the following three lemmas.

\medskip

\begin{lem}\label{many_non-commensurable} {\rm(see~\cite[Lemma~10.3]{Bog_1})}
Let $G$ be a group, $X\subseteq G$, $H\hookrightarrow_h(G,X)$ a finitely generated infinite subgroup.
Then for any finite collection of elements $a_1,\dots,a_s\in G\setminus H$ and any infinite
subset $\widetilde{H}\subseteq H$, there exist elements $h_1,\dots ,h_s\in \widetilde{H}$
such that $a_1h_1,\dots, a_sh_s$ are paarwise non-commensurable
loxodromic elements with respect to the action of $G$ on $\Gamma(G,X\sqcup H)$.
\end{lem}

\medskip

\noindent
{\it Remark.} We will use a special case of this lemma where $a_1=\dots=a_s$.
In this case this lemma is very similar to~\cite[Corollary 6.12]{DOG} and, actually, can be deduced
from the proof of this corollary.

%\medskip

\begin{lem}\label{many_special} {\rm(\cite[Lemma 10.4]{Bog_1})}
Let $G$ be an acylindrically hyperbolic group with respect to a generating set $Y$ and let $a,b\in G$ be
two non-commensurable loxodromic with respect to $Y$ elements, where, additionally, $a$ is special.
Then there exists a positive integer $n_0$ such that for any $n,m\geqslant n_0$ the element $g=a^nb^m$
is special with respect to some generating set, in particular $E_G(g)=\langle g\rangle$.
\end{lem}

\begin{lem}\label{one_special}{\rm(\cite[Lemma 10.5]{Bog_1})}
Suppose that $G$ is an acylindrically hyperbolic group
without nontrivial finite normal subgroups. Then
%$G$ contains at least one special element. Moreover,
there exist an element $g\in G$ and a generating set $Y$ of $G$ such that
$g$ is special with respect to $Y$ and $\langle g\rangle\hookrightarrow_h(G,Y)$.
\end{lem}

\begin{prop}\label{many_jointly_special}
Suppose that $G$ is an acylindrically hyperbolic group without nontrivial finite normal subgroups.
Then, there are special loxodromic elements $a,g\in G$ such that for any integer $k>0$,
the coset $a\langle g\rangle$ contains $k$ pairwise non-commensurable and jointly special elements.
\end{prop}

\medskip

{\it Proof.} By Lemma~\ref{one_special},
there exist an element $g\in G$ and a generating set~$Y$ of $G$ such that
$g$ is special with respect to $Y$
and $H\hookrightarrow_h(G,Y)$, where $H=\langle g\rangle$.
It follows from Definition~\ref{exact_defn_special} that $G$ is acylindrically hyperbolic with respect to $Y$.

Let $b\in G\setminus H$ be an arbitrary element. By Lemma~\ref{many_non-commensurable},
there exist two non-commensurable loxodromic elements $bg^s$, $bg^t$ with respect to $Y\sqcup H$.\break
It follows that they are loxodromic with respect to $Y$.
At least one of them, say $c:=bg^s$ is non-commensurable with $g$.
In particular,
$$
\langle c\rangle \cap \langle g\rangle=1.
$$
By Lemma~\ref{many_special},
there exists a positive integer $n_0$ such that for any $n,m\geqslant n_0$ the element $c^ng^m$
is special with respect to some generating set, in particular
$$
E_G(c^ng^m)=\langle c^ng^m\rangle.\eqno{(4.1)}
$$

By Lemma~\ref{many_non-commensurable}, for any $k\in \mathbb{N}$, there exist natural numbers $m_1<m_2<\dots <m_{k+1}$ such that $m_1\geqslant n_0$ and the elements
$c^{n_0}g^{m_1},\dots,c^{n_0}g^{m_{k+1}}$ are pairwise non-commensurable and loxodromic with respect to $Y\sqcup H$. Then they are loxodromic with respect to $Y$.
Moreover, by (4.1), we have $E_G(c^{n_0}g^{m_i})=\langle c^{n_0}g^{m_i}\rangle$ for $i=1,\dots,k+1$.
Thus, these elements are pairwise non-commensurable and jointly special with respect to $Y$.

We set $a=c^{n_0}g^{m_1}$. Then $a$, $g$ and $k$ elements
$ag^{m_2-m_1},\dots ,ag^{m_{k+1}-m_1}$
%the set $I=\{m_2-m_1,\dots,m_{k+1}-m_1\}$
satisfy the conclusion of proposition.\hfill $\Box$

\medskip

{\it Remark.} One can prove a stronger version of this lemma, saying that for any infinite subset $I\subseteq a\langle g\rangle$, there exists an infinite subset of $I$ consisting of pairwise
non-commensurable and jointly special elements.

\section{Test words in acylindrically hyperbolic groups}

The history of test words in free groups is illuminated in~\cite{Ivanov}. In this paper Ivanov constructed the so-called $C$-test words in free groups. In~\cite{Lee} Lee constructed $C$-test words with some additional property.
%In~\cite{MR}, Myasnikov and Roman'kov used Lee's test words to study verbally closed subgroups of free groups.
In this section we construct certain test words in acylindrically hyperbolic groups.

\begin{defn}\label{def_text-word} {\rm (see~\cite[Definition 12.1]{Bog_1})}
{\rm Let $H$ be a group and let $a_1,\dots,a_k$ be some elements of~$H$.
A word $W(x_1,\dots ,x_k)$ from $F_k$ is called an
{\it $(a_1,\dots,a_k)$-test word} if for every solution $(b_1,\dots,b_k)$ of the equation
$$
W(a_1,\dots ,a_k)=W(x_1,\dots ,x_k)
$$
in $H$, there exists a number $\alpha\in \mathbb{Z}$ such that
$b_i=a_i^{U^{\alpha}}$ for $i=1,\dots, k$, where $U=W(a_1,\dots ,a_k)$.
}
\end{defn}

{\bf Notation.} We write $\bold{1}^k$ for the tuple $(\underbrace{1,\dots,1}_{k})$.

\medskip

The following proposition is a special case of Proposition 12.4 from~\cite{Bog_1} .

\begin{prop}\label{test_special_words} {\rm (see~\cite[Proposition 12.4]{Bog_1})}
Let $H$ be an acylindrically hyperbolic group without nontrivial finite normal subgroups
and let $a_1,\dots, a_k\in H$ (where $k\geqslant 3$) be jointly special and pairwise non-commensurable elements.
Then there is an $(a_1,\dots,a_k,\bold{1}^{k-2})$-test word $W_k(x_1,\dots,x_k,y_3,\dots,y_k)$.

Moreover, one can choose this test word so that the elements $a_1,\dots,a_k$ together with
$W_k(a_1,\dots,a_k,\bold{1}^{k-2})$ are jointly special and pairwise non-commen\-surable.
\end{prop}

%The first part of this proposition is used in the proof of the following theorem.

\medskip

In the following section we need a weaker version of Proposition~\ref{test_special_words},
which will also simplify notations.
We set
$$
\mathcal{U}_k(x_1,\dots,x_k):=W_k(x_1,\dots,x_k,\bold{1}^{k-2}).
$$
Then we obtain the following corollary.

\begin{cor}\label{test_special_words_1}
Let $H$ be an acylindrically hyperbolic group without nontrivial finite normal subgroups
and let $a_1,\dots, a_k\in H$ (where $k\geqslant 3$) be jointly special and pairwise non-commensurable elements.
Then there is an $(a_1,\dots,a_k)$-test word $\mathcal{U}_k(x_1,\dots,x_k)$ such that
the elements $a_1,\dots,a_k$ together with
$\mathcal{U}_k(a_1,\dots,a_k)$ are jointly special and pairwise non-commensurable.
\end{cor}

\medskip

\section{Proof of Theorem A}

%\subsection{Coding a finite system of equations by one equation}

\begin{lem}\label{splitted_equations}
Let $H$ be a group. For any finite system of equations $S\subset F_n\ast H$
%By introducing new variables and equations of the type $z_i=c_i$, where $c_i\in H$
there exists a finite system $S'\subset F_n\ast H$ consisting of only splitted equations
such that $|S'|\geqslant |S|$ and
%the elements $h_1,h_2$ are non-commensurable and jointly special, and
$$
V_H(S')=V_H(S)\times \{g_1\}\times \dots \times \{g_k\}
$$
for some elements $g_1,\dots,g_k\in H$.
\end{lem}

\medskip

{\it Proof.} To define $S'$, one should replace each constant $h$ in $S$ by a new variable $x_h$ and add the equation $x_hh^{-1}$. \hfill $\Box$

\bigskip

\noindent
{\bf Notation.} To shorten notation, we write $\underline{x}$ instead of the tuple $(x_1,\dots,x_n)$.

\bigskip

{\it Proof of Theorem A.}

(1) Let $S=\{s_1,\dots ,s_k\}$, where $s_i\in F_n\ast H$.
We take arbitrary $(k+2)$ jointly special and pairwise non-commensurable elements $a_1,\dots, a_{k+2}\in H$.
The existence of such elements is guaranteed by Proposition~\ref{many_jointly_special}.
By Proposition~\ref{test_special_words}, there exists an $(a_1,\dots ,a_{k+2},\bold{1}^{k})$-test word
$W_{k+2}(z_1,\dots,z_{k+2},y_1,\dots,y_k)$.
Then the desired equation is
$$f:\hspace*{2mm} W_{k+2}(a_1,\dots ,a_{k+2},\bold{1}^{k})=W_{k+2}(a_1,\dots ,a_{k+2},s_1,\dots,s_k).$$

\medskip

(2) By Lemma~\ref{splitted_equations}, we may assume that
$S$ consists of splitted equations:
$$
S=\{w_i(x_1,\dots ,x_n)=h_i\,|\, i=1,\dots,m\}.\eqno{(6.1)}
$$
By Proposition~\ref{many_jointly_special}, there exist special elements $a,b$ such that
the coset $a\langle b\rangle$ contains $(2m+2)$ pairwise non-commensurable jointly special elements,
say
$$
 ab^s,ab^t,\, ab^{k_1},\dots,ab^{k_m},ab^{l_1},\dots,ab^{l_m}.
$$
Then the elements of the tuple
$$
T=(ab^s,ab^t,\,\,h_i^{-1}ab^{k_i}h_i,\,\,h_i^{-1}ab^{l_i}h_i;\,\, i=1,\dots,m)
$$
are also pairwise non-commensurable and jointly special.
Let $\mathcal{U}_{2m+2}$ be the $T$-test word, see Corollary~\ref{test_special_words_1}. We set

$$
f_0=
\mathcal{U}_{2m+2}(ab^s,ab^t,\,\,h_i^{-1}ab^{k_i}h_i,\,\,h_i^{-1}ab^{l_i}h_i;\,\, i=1,\dots,m).
$$

\medskip

Now we introduce $2$ new variables $y,z$ and set

$$
f_1^{-1}=
\mathcal{U}_{2m+2}(yz^s,yz^t,\,\,w_i(\underline{x})^{-1}yz^{k_i}w_i(\underline{x}),\,\, w_i(\underline{x})^{-1}yz^{l_i}w_i(\underline{x});\,\, i=1,\dots,m).
$$

\medskip
\noindent
We show that the splitted equation $f$ written in the form
$$
f_1^{-1}=f_0.\eqno{(6.2)}
$$
satisfies the statements (a) and (b) of Theorem {\bf A}.

\medskip

(a)
Suppose that $f$ has a solution in $H$, say
$$
(x_1,\dots,x_n,y,z)=(C_1,\dots,C_n,A,B).
$$

We shall show that there exists $\alpha\in \mathbb{Z}$ such that $(C_1,\dots,C_n)^{f_0^{-\alpha}}$
is a solution of the system $S$.

Set $H_i=w_i(\underline{C})$.
Then we have
$$
\begin{array}{ll}
& \mathcal{U}_{2m+2}(ab^s,ab^t,\,\,h_i^{-1}ab^{k_i}h_i,\,\,h_i^{-1}ab^{l_i}h_i;\,\, i=1,\dots,m)\\
= & \\
& \mathcal{U}_{2m+2}(AB^s,AB^t,\,\,H_i^{-1}AB^{k_i}H_i,\,\,H_i^{-1}AB^{l_i}H_i;\,\, i=1,\dots,m).
\end{array}
\eqno{(6.3)}
$$

\medskip

By definition of the test word (see Definition~\ref{def_text-word}), applied to $\mathcal{U}_{2m+2}$ and (6.3), there exists $\alpha\in \mathbb{Z}$ such that the formulas (6.4) -- (6.7) are valid:

$$
(ab^s)^{f_0^{\alpha}}=AB^s,\eqno{(6.4)}
$$

$$
(ab^t)^{f_0^{\alpha}}=AB^t,\eqno{(6.5)}
$$

$$
(h_i^{-1}ab^{k_i}h_i)^{f_0^{\alpha}}=H_i^{-1}AB^{k_i}H_i,\eqno{(6.6)}
$$

$$
(h_i^{-1}ab^{l_i}h_i)^{f_0^{\alpha}}=H_i^{-1}AB^{l_i}H_i,\eqno{(6.7)}
$$

\medskip

It follows from (6.4) and (6.5) that $(b^{s-t})^{f_0^{\alpha}}=B^{s-t}$.
Since $b$ is special, we have
$$
B=b^{f_0^{\alpha}}.\eqno{(6.8)}
$$

\medskip

From (6.4) and (6.8), we obtain that
$$
A=a^{f_0^{\alpha}}.\eqno{(6.9)}
$$

Substituting (6.8) and (6.9) in (6.6) and (6.7), we deduce that
%\medskip
%\medskip

%It follows from (5.12) and (5.13) that
$$
\begin{array}{ll}
h_i^{f_0^{\alpha}}H_i^{-1} & \in C_H((ab^{k_i})^{f_0^{\alpha}})\cap C_H((ab^{l_i})^{f_0^{\alpha}})\vspace*{3mm}\\
& \overset{(3.1)}{\subseteq} E_H((ab^{k_i})^{f_0^{\alpha}})\cap E_H((ab^{l_i})^{f_0^{\alpha}})\vspace*{3mm}\\
& = \langle (ab^{k_i})^{f_0^{\alpha}}\rangle \cap \langle(ab^{l_i})^{f_0^{\alpha}}\rangle
\end{array}
$$
The latter equation holds since the elements $ab^{k_i}$ and $ab^{l_i}$ are special.
This intersection is trivial since $ab^{k_i}$ and $ab^{l_i}$ are non-commensurable.
Therefore
$$
w_i(\underline{C})=h_i^{f_0^{\alpha}}
$$
for $i=1,\dots,m$. Hence $\underline{C}^{f_0^{-\alpha}}$ is a solution of $S$.

\medskip

(b)
Since $V_G(f)$ is invariant under conjugation by the element $f_0$, it suffices to check that
$$
{\text{\bf pr}}_n\bigl(V_G(f)\bigr)\supseteq V_G(S).
$$

The latter is trivial: If $(x_1,\dots,x_n)=(c_1,\dots,c_n)$ is a solution of the system $S$ in $G$,
then
$$
(x_1,\dots,x_n,y,z)=(c_1,\dots,c_n,a,b)
$$
is a solution of the equation (6.2) in $G$.
%\hfill $\Box$

%\subsection{Coding a finite system of equation by two splitted equations}

%\begin{thm}\label{algebraic_sets_1}
%Let $H$ be an acylindrically hyperbolic group without nontrivial finite normal subgroups.
%For any finite system of equations $S\subseteq F_n\ast H$ there exist two splitted equation
%$f,g\in F_{6\ell(S)+2}\ast H$ such that
%$$
%V_H(S)={\text{\bf pr}}_n\bigl(V_H(f)\bigr)\cap {\text{\bf pr}}_n\bigl(V_H(g)\bigr).
%$$
%\end{thm}

\medskip

%{\it Proof.}

(3)
We may again assume that $S$ has the form (6.1).
We may additionally assume that the set $\{h_1,\dots,h_m\}$ of right sides of equations from $S$
contains two non-commensurable special elements from $H$; otherwise we could take two non-commensurable special elements $u,v\in H$ and add two equations $x_{n+1}=u$ and $x_{n+2}=v$ to $S$.
Obviously, the set of solutions of the old system~$S$ is a projection of the set of solutions of the new system $S$.

In the following we will use the tuple $T$, the element $f_0$ and the equation $f$ defined in (a).
Thus, we have
$$
{\text{\bf pr}}_n\bigl(V_H(f)\bigr)=\underset{\alpha\in \mathbb{Z}}{\bigcup} V_H(S)^{f_0^{\alpha}}.
\eqno{(6.10)}
$$

By Corollary~\ref{test_special_words_1}, all components of $T$ together with $f_0$
are pairwise non-commensurable and jointly special.
Let $T'$ be the tuple obtained from $T$ by adding the component~$f_0$:
$$
T'=(f_0,\,ab^s,ab^t,\,\,h_i^{-1}ab^{k_i}h_i,\,\,h_i^{-1}ab^{l_i}h_i;\,\, i=1,\dots,m)
$$

Let $\mathcal{U}_{2m+3}$ be the $T'$-test word from Corollary~\ref{test_special_words}.
We set

$$
g_0=\mathcal{U}_{2m+3}(f_0;\,
ab^s,ab^t,\,\,h_i^{-1}ab^{k_i}h_i,\,\,h_i^{-1}ab^{l_i}h_i;\,\, i=1,\dots,m).
$$

\medskip

Now we introduce new variables $t,y,z$ and define the word

$$
g_1^{-1}=\mathcal{U}_{2m+3}(t,\, yz^s,yz^t,\,\,w_i(\underline{x})^{-1}yz^{k_i}w_i(\underline{x}),\,\, w_i(\underline{x})^{-1}yz^{l_i}w_i(\underline{x});\,\, i=1,\dots,m).
$$
Let $g$ be the equation $g_1g_0$. Using the same arguments as in the proof of (a), we obtain
$$
{\text{\bf pr}}_n\bigl(V_H(g)\bigr)=\underset{\alpha\in \mathbb{Z}}{\bigcup} V_H(S)^{g_0^{\alpha}}.
\eqno{(6.11)}
$$

\medskip

{\bf Claim 1.} Suppose that for some $\alpha,\beta\in \mathbb{Z}$ we have
$$
V_H(S)^{f_0^{\alpha}}\cap V_H(S)^{g_0^{\beta}}\neq \emptyset.
$$
Then $\alpha=\beta=0$.

\medskip

{\it Proof.}
By assumption, the set $\{h_1,\dots,h_m\}$ of right sides of equations from $S$ contains two non-commensurable special elements,
say $u,v$. Then
$$
u^{f_0^{\alpha}}=u^{g_0^{\beta}}\hspace*{2mm} {\text{\rm and}}\hspace*{2mm}v^{f_0^{\alpha}}=v^{g_0^{\beta}},
$$
and we deduce $$f_0^{\alpha}g_0^{-\beta}\in E_H(u)\cap E_H(v)=\langle u\rangle\cap \langle v\rangle=1.$$
The penultimate equation holds since $u$ and $b$ are special, and the latter equation holds since $u$ and $v$
are non-commensurable.

By Corollary~\ref{test_special_words_1}, all components of $T'$ together with $g_0$ are pairwise non-commen\-surable and jointly special.
In particular, $g_0$ and $f_0$ are non-commen\-surable and have infinite orders.
From this and $f_0^{\alpha}=g_0^{\beta}$, we obtain $\alpha=\beta=0$.
\hfill $\Box$

\medskip

This claim and equations (6.10) and (6.11) imply that
$$
V_H(S)={\text{\bf pr}}_n\bigl(V_H(f)\bigr)\cap {\text{\bf pr}}_n\bigl(V_H(g)\bigr).
$$

\hfill $\Box$

\begin{rmk}
{\rm
In the general case, one splitted equation in statement (3) of Theorem {\bf A} is not sufficient.
Indeed,
let $H$ be an arbitrary nontrivial group and let
$
S=\{w_i(x_1,\dots ,x_n)=h_i\,|\, i=1,\dots,m\},
$
$m\geqslant 2$, be a finite system of splitted equations with constants $h_i$ from $H$  such that
$C_H(h_1)\cap C_H(h_2)=1$ and $V_H(S)\neq \emptyset$. Then for any splitted equation $f\in F_k\ast H$, where $k\geqslant n$, we have
$$
V_H(S)\neq{\text{\bf pr}}_n\bigl(V_H(f)\bigr).
$$
This follows from the following obsevations:
%statement (3) of Theorem B, we use two splitted equations.
%We explain, why one splitted equation is not sufficient in general.
\begin{enumerate}
\item[(a)] If $f$ is a splitted equation of the form $f_1f_0$, where $f_1\in F_k$ and $f_0\in H$, then
$\bigl(V_H(f)\bigr)^{f_0}=V_H(f)$.
Moreover, we have $\bigl(V_H(f)\bigr)^{g}=V_H(f)$ for any $g\in H$ if $f_0=1$.

\medskip

\item[(b)] %If $S$ is a system of splitted equations containing two non-commensurable special elements, then
$V_H(S)^g\cap V_H(S)=\emptyset$ for every nontrivial $g\in H$.
This can be proved similarly to the proof of Claim 1.
\end{enumerate}
}
\end{rmk}

\section{Proof of Theorem C}

\medskip

%\noindent
%{\bf Theorem C.}
%{\it Let $H$ be an acylindrically hyperbolic group without nontrivial finite normal subgroups
%and let $G$ be an arbitrary overgroup of $H$.
%Then $H$ is verbally closed in $G$ if and only if $H$ is algebraically closed in $G$.
%}

%\medskip

{\it Proof.}
Suppose that $H$ is verbally closed in $G$. We show that $H$ is algebraically closed in $G$.
%We will prove that $H$ is algebraically closed in $G$.
Let $S$ be a finite system of equations with constants from $H$ such that $V_G(S)\neq  \emptyset$.
We shall show that $V_H(S)\neq  \emptyset$.

Let $f$ be a splitted equation as in statement (2) of Theorem {\bf A}.
By part (b) of this statement, we have $V_G(S)\subseteq {\text{\bf pr}}_n \bigl(V_G(f)\bigr)$, hence
$V_G(f)\neq  \emptyset$. Since $H$ is verbally closed in $G$, we have $V_H(f)\neq  \emptyset$.
By part (a) of statement (2) of Theorem {\bf A}, we have $V_H(S)\neq  \emptyset$. Thus, $H$ is algebraically closed in $G$.
The converse implication is obvious.
\hfill $\Box$

\begin{rmk}\label{Remark_free_prod}
{\rm
%1)
%Let $H$ be a virtually free group without nontrivial finite normal subgroups and let $G$ be an overgroup of $G$.
%In~\cite{KMM}, Klyachko, Mazhuga and Miroshnichenko proved the following result:
%Let $H$ be a virtually free group without nontrivial finite normal subgroups and
%let $G$ be an overgroup of $G$. If $H$ is verbally
%The main part of their paper is devoted to considering the case $H=\mathbb{Z}_2\ast \mathbb{Z}_2$.
Consider the free product $$H=\underset{\alpha\in \frak{A}}{\ast} H_{\alpha},$$ where $\frak{A}$ is an arbitrary set of cardinal larger than 1 and each $H_{\alpha}$ is nontrivial.
In~\cite{Mazhuga_3} Mazhuga showed that if $H$ is a verbally closed subgroup of a group $G$,
then $H$ is algebraically closed in $G$.
This result (except the very special case $H=\mathbb{Z}_2\ast \mathbb{Z}_2$, which was first considered in~\cite{KMM}) follows from our Theorem {\bf C}.

Indeed, $H$ can be splitted as $H=A\ast B$, where $A$ and $B$ are nontrivial;
hence $H$ is relatively hyperbolic with respect to $\{A, B\}$.
%The group $H$ is virtually cyclic if and only if it has the form $H\cong\mathbb{Z}_2\ast \mathbb{Z}_2$.
It is well known that if a non-(virtually cyclic) group is relatively hyperbolic with respect to a collection of proper subgroups, then it is acylindrically hyperbolic.
Therefore $H$ is acylindricaly hyperbolic except the case, where
$H\cong\mathbb{Z}_2\ast \mathbb{Z}_2$.
%Note that this case was first considered in~\cite{KMM}.
}
\end{rmk}

%/////////////////////

%\begin{cor}
%Let $H$ be a finitely generated acylindrically hyperbolic group without nontrivial finite normal subgroups.
%Then there exists an element $h\in H$ such that if $\phi$ an $\psi$ two endomorphism of $H$ satisfying
%$\phi(h)=\psi(h)$, then $\phi=\widehat{g}\circ \psi$ for some $a\in H$.
%\end{cor}

%\begin{cor}
%Let $G$  be a finitely generated group.
%Let $H$ be a finitely generated acylindrically hyperbolic group without nontrivial finite normal subgroups.
%Then there exists an element $g\in G$ with the following property:
%If $\phi,\psi\in Hom(G,H)$, then $\phi\sim \psi$ iff $\phi(g)$ is conjugate to $\psi(g)$.
%\end{cor}

\end{document}